\documentclass[
  a4paper, 
  reqno, 
  oneside, 
  11pt
]{amsart}

\usepackage{amsmath,amssymb,bm,cancel,booktabs, epsfig}
\usepackage{hyperref}
\usepackage[dvipsnames]{xcolor}
\usepackage{tikz,bbding}  % for diagrams
\usetikzlibrary{calc,decorations.pathmorphing}
\usetikzlibrary{arrows,positioning}
\usepackage{graphicx,caption}

\theoremstyle{plain}	% 'plain' is the default.  The others are 'definition' and 'remark'.
\newtheorem*{theorem*}{Theorem}

\newtheorem*{conjecture*}{Conjecture}
\theoremstyle{remark}

\DeclareMathOperator{\Tot}{Tot}
\DeclareMathOperator{\Ad}{Ad}

\DeclareMathOperator{\Sing}{Sing}

\newcommand{\ce}{\mathrel{\mathop:}=}

\newcommand{\bcomment}[1]{{\color{olive} \footnotesize #1}}

\title[Curvature Grafted by  Instantons]{Curvature Grafted  by  Instantons}
\author{Elizabeth Gasparim {\tiny and} Bruno Suzuki}
\address{Depto. Matem\'aticas, Universidad Cat\'olica del Norte, Chile, etgasparim@gmail.com, obrunosuzuki@gmail.com}
\thanks{B.S. acknowledges  support of 
ANID 2019/13204-0, E.G.  acknowledges  support of VRIDT-UCN, Chile}

\newcommand{\trevo}{\textnormal{\tiny{\FourClowerSolid}}}
\newcommand{\atrevo}{\tiny{\textnormal{\FourClowerOpen}}}

\begin{document} 

\maketitle

\centerline{\emph {In memory of  Pushan Majumdar.}}

\begin{abstract} We show that an instanton with high charge can provoke the 
creation of extra curvature on the space that holds it. Geometrically, this corresponds 
to a new surgery operation, which we name grafting. Curvature around a sphere increases
by grafting when the charge of an instanton decays. 
\end{abstract}

\section{Particle decay and local characteristic classes}

The concept of local characteristic classes has its motivation in the description of 
particle decay phenomena. 
 Recent news on particle decay 
have brought to light the need to describe a general theory of local 
invariants. 
The new observation of decay  was presented at the 40th International Conference on High Energy Physics. 
Scientists at CERN have reported on their first significant evidence for a process predicted by theory, paving the way for searches for evidence of new physics in particle processes that could explain dark matter and other mysteries of the universe.
 CERN NA62 presented the first significant experimental evidence for the ultra-rare decay of 
 the charged kaon into a charged pion and two neutrinos.

\begin{center}
\begin{tikzpicture}
\newcommand{\rr}{.9}

\coordinate (O) at (0,0);
\coordinate (C1) at (-0.5,0);
\coordinate (C2) at (0.5,0);
\coordinate (V) at (-2,-3);
\coordinate (V2) at (2,-3);

%kaon
{
	\draw[fill, gray] (O) circle [radius=\rr+.2];

	\draw[fill, yellow] (C1) circle (\rr/3);
	\draw[fill, Cerulean] (C2) circle (\rr/3);

	\draw[YellowOrange, thick, decorate, decoration={snake, segment length=1.5mm, amplitude=2mm}] ($(C1)+(.2,0)$) -- ($(C2)-(.2,0)$); %snakeline
	
	\node at (C1) {$\overline{u}$};
	\node at (C2) {$s$};
}

%pi­on
{
	\coordinate (O2) at ($(O)+(V)$);
	\coordinate (D1) at ($(C1)+(V)$);
	\coordinate (D2) at ($(C2)+(V)$);

	\draw[fill, gray] (O2) circle [radius=\rr];

	\draw[fill, Cerulean] (D1) circle (\rr/3);
	\draw[fill, yellow] (D2) circle (\rr/3);

\draw[YellowOrange, thick, decorate, decoration={snake, segment length=1.5mm, amplitude=2mm}] ($(D1)+(.2,0)$) -- ($(D2)-(.2,0)$); %snakeline

	\node at (D1) {$u$};
	\node at (D2) {$\overline{d}$};
}

%neutrino1
{
	\coordinate (N1) at ($(O)+(0,-3)$);
	\draw[fill, Tan] (N1) circle [radius=\rr/3];
	\node at (N1) {$n$};
}

%neutrino2
{
	\coordinate (N2) at ($(O)+(V2)$);
	\draw[fill, Tan] (N2) circle [radius=\rr/3];
	\node at (N2) {$n$};
}

%arrows
{
	\draw [->, thick] ($0.35*(V)$) -- ($0.7*(V)$);
	\draw [->, thick] ($0.4*(N1)$) -- ($0.7*(N1)$);
	\draw [->, thick] ($0.35*(N2)$) -- ($0.7*(N2)$);
}

\end{tikzpicture}
\end{center}

 Here we will  discuss the concepts of local characteristic classes, in particular Chern class 
and local holomorphic Euler characteristic, and how these occur within the  mathematical descriptions of particle decay.
 We will focus in 4 dimensions, recalling results obtained in joint work with Pushan Majumdar about instanton decay, and we will propose generalizations to higher dimensions.
 We are especially interested in the loss of charge that can be provoked by the contraction of a 2-sphere. The local Chern classes measure local manifestations of curvature concentrated around points or around subvarieties.

\begin{center}
\includegraphics[height=4.2cm]{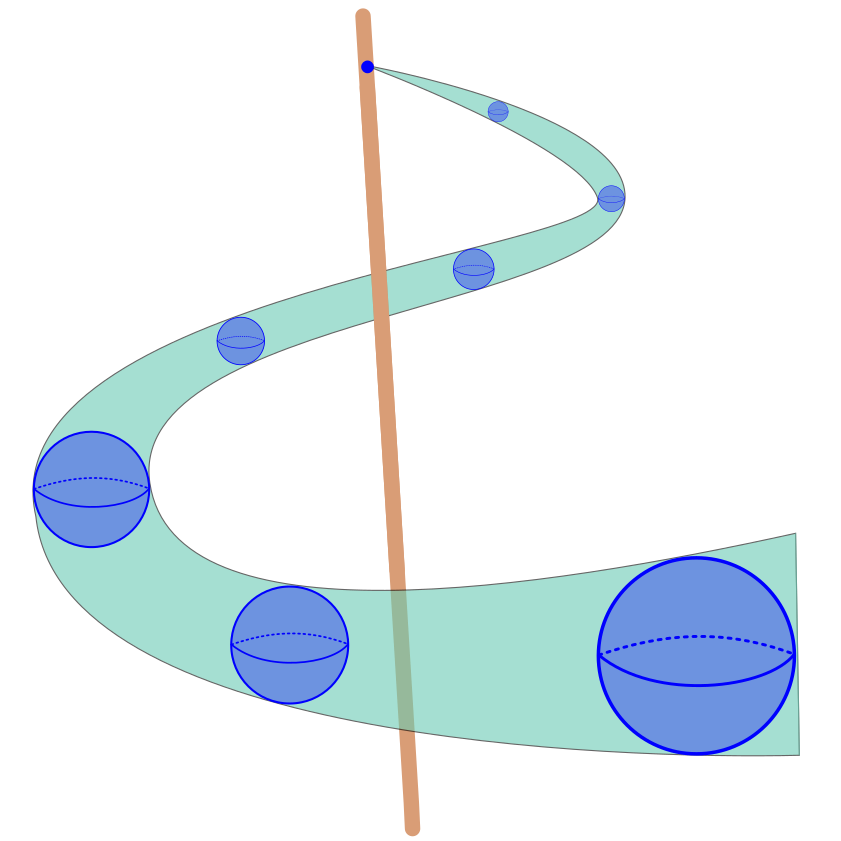}
\end{center}

\section{Instantons}
We first review the basic mathematical setup of instantons. 
The original  definition of  instanton says it is  a connection minimizing the Yang--Mills 
functional. Let us make this concept mathematically precise. 

We start up by recalling the concept of connection on a principal bundle.
Let  $P$   be a  principal $G$-bundle over a smooth manifold $M$.
 A \emph{connection} on $P$ is a differential 1-form $\omega$ on $P$ 
 with values in the Lie algebra satisfying:
 \begin{itemize}
 \item $\Ad_gR_g^*\omega = \omega$ ($G$-equivariance),
 \item $\omega(X_\xi)=\xi$ (Fundamental vector fields),
 \end{itemize}
where $R_g$ denotes the right translation by $g$, upper $*$ the pull back of forms,
and $X_\xi$ denotes the fundamental vector field associated to  $\xi$.

Let $F_A$ denote the curvature of the connection $A$, and assume $M$ has real dimension $4$.
An \emph{instanton}  on $M$ is a  connection minimizing the \emph{Yang--Mills Functional}:
$$YM(A) = \int |F_A|^2 \textnormal{dvol}.$$
Note that the functional applies to 4 dimensions, since the integrand is a 4-form.
The \emph{Yang--Mills equations}
 are the Euler--Lagrange equations corresponding to the Yang--Mills Functional.
Hence, an instanton is a solution of the YM equations.
These  are nonlinear partial differential equations, which are well known to be very difficult to solve.

One of the few cases where an explicit solution is written in terms of connections is the 
basic $SU(2)$ instanton on $S^4$:
$$A= \frac{1}{1+|x|^2}(\theta_1{\bf i}+ \theta_2{\bf j}+\theta_3{\bf k})$$
where $ {\bf i}, {\bf j}, {\bf k}$ is a basis for $\mathfrak{sl}(2)$ 
and 
$$\theta_1= x_1\textrm{d}x_2 -x_2\textrm{d}x_1-x_3\textrm{d}x_4+x_4\textrm{d}x_3,$$ 
$$\theta_2= x_1\textrm{d}x_3 -x_3\textrm{d}x_1-x_4\textrm{d}x_2+x_2\textrm{d}x_4,$$ 
$$\theta_3= x_1\textrm{d}x_4 -x_4\textrm{d}x_1-x_2\textrm{d}x_3+x_3\textrm{d}x_2.$$ 
\bcomment{}
However, for more general manifolds explicit connections solving the YM equations are not known. 
 A simpler way to look for instantons is to consider  the linearized version of the YM equations, namely:
$$\star F = -F,$$ where $\star$ is the Hodge operator. 
Connections whose curvatures satisfy this equation are called \emph{anti-self-dual}.
 These ASD connections are solutions of the YM equations. 
 Using ASD connections, alternative approaches 
to find instantons appeared.  
 A celebrated construction of instantons in Euclidean space was obtained by Atiyah--Drinfeld--Hitchin--Manin
in \cite{ADHM} in terms of the following linear data:
\begin{itemize}
\item complex vector spaces $V$ and $W$ of dimension $k$ and $N$,
\item $k \times k$ complex matrices $B_1$, $B_2$, 
\item a $k \times N$ complex matrix $I$, 
\item a $N \times k$ complex matrix $J$,
\item a real moment map 
${\displaystyle \mu _{r}=[B_{1},B_{1}^{\dagger }]+[B_{2},B_{2}^{\dagger }]+II^{\dagger }-J^{\dagger }J,}$
\item a complex moment map 
${ \displaystyle \mu _{c}=[B_{1},B_{2}]+IJ.}$
\end{itemize}
The ADHM construction was explored and generalized in a variety of ways by several authors. 
Subsequently,  this approach using complex vector spaces was significantly improved  as to 
include  a far wider variety of 
instantons  by using  holomorphic vector bundles, giving rise to  
a very general principle called the 
\emph{Kobayashi--Hitchin correspondence}, see \cite{LT}, which provides an identification: 
\begin{center} ASD $\mathrm{SU}(r)$  connections $\leftrightarrow$ rank $r$ holomorphic bundles.
\end{center}
The correspondence was proved by Donaldson \cite{D} and Uhlenbeck--Yau \cite{UY}.
We will concern ourselves with the case of $\mathrm{SU}(2)$ instantons:

\[
\left\{
\begin{array}{c}
\text{irreducible $\mathrm{SU} (2)$-instantons} \\ 
\text{of charge $n$}
\end{array}
\right\}
\leftrightarrow
\left\{
\begin{array}{c}
\text{stable $\operatorname{SL} (2, \mathbb C)$-bundles} \\
\text{with $c_2 = n$}
\end{array}
\right\}
\]\\

If we represent the connection in terms of a covariant derivative on a vector bundle, then 
the association is:
$$\nabla= \partial + \bar\partial \phantom{xx}\leftrightarrow  \phantom{xxx}\bar\partial.$$

The method of using holomorphic bundles to study instantons is very fruitful, and has 
been used widely in the study of instanton moduli spaces, as it allows us to profit from all the 
technology of complex algebraic geometry to explore properties of instantons. 
In what follows we will identify instantons with holomorphic bundles via the 
Kobayashi--Hitchin  correspondence, using the particular version of this conjecture 
proved in \cite{GKM} for instantons on the surfaces $Z_k$.
  We can then talk about instanton charges in terms of characteristic classes.

\section{Chern classes}

Let us recall the basic properties of Chern classes of vector bundles.
Let $E$ be a vector bundle on a complex manifold $M$.  The \emph{Chern polynomial} of $E$ is
written as
$$c(E)= c_0(E)+c_1(E)t+ c_2(E)t^2 + \cdots + c_n(E)t^n$$
where $c_i(E) \in H^{2i}(M,\mathbb Z)$ is the $i$-th Chern class of $E$.
These Chern classes satisfy the following 4 axioms:
\begin{enumerate}
\item $c_0(E)=1$
\item $c_i(f^*(E))= f^*(c_i(E))$
\item $ 0 \rightarrow F_1 \rightarrow E \rightarrow F_2\rightarrow 0 \Rightarrow c(E) = c(F_1)c(F_2)$
\item $c_1(\gamma) =-1,$
\end{enumerate}
where $\gamma $ is the universal bundle on $\mathbb P^1$ whose fiber at 
a given point is the line in $\mathbb C^2$ it represents. 

Observe that  the single axiom setting the first Chern class of $\gamma$ to $-1$ is the only one among the 4
axioms that guaranties nontriviality of the theory, else, without this axiom, the other properties would all have been 
satisfied by $c=1$ for all bundles. The line bundle $\gamma$ is usually denoted in algebraic 
geometry by $ \mathcal O_{\mathbb P^1}(-1)$ and its total space is isomorphic to 
the blow-up of $\mathbb C^2$ at the origin: $\Tot \mathcal O_{\mathbb P^1}(-1) \simeq \widetilde{\mathbb C^2}$, 
where the blow-up of the complex plane is:
$$\widetilde{\mathbb C^2}\ce \{(z,l) \in \mathbb C^2 \times \mathbb P^1 : z \in l\}.$$

We will return to $\widetilde{\mathbb C^2}$ in the section about local surfaces, where it will be denoted by 
$Z_1$ being the first in the collection of surfaces $Z_k\ce \Tot \mathcal O_{\mathbb P^1}(-k)$ 
we will use, where 
$\mathcal O_{\mathbb P^1}(-k) = \mathcal O_{\mathbb P^1}(-1)\otimes 
\cdots \otimes \mathcal O_{\mathbb P^1}(-1)$  tensored $k$ times has first Chern class 
$c_1(\mathcal O_{\mathbb P^1}(-k)) = -k$ by  axiom 3. 
Given its importance in the context of characteristic classes and also for 
our calculations of local Chern classes, it is worth depicting the geometry of $Z_1$ 
in terms of the real analogue. The blow up of the real plane is 
$$ \widetilde{\mathbb R^2}\ce \{(z,l) \in \mathbb R^2 \times \mathbb R\mathbb P^1 : z \in l\}$$
which can be depicted as
\begin{center}\includegraphics[height=6.2cm]{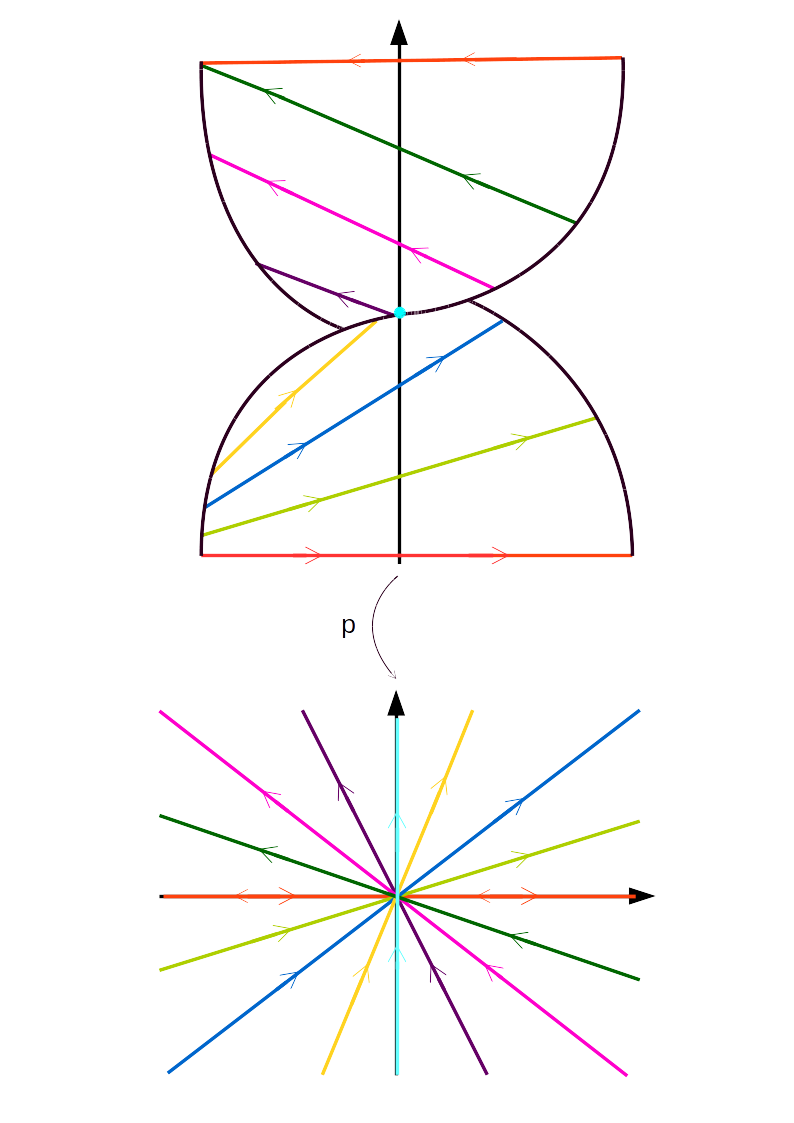}\end{center}
and is isomorphic to an infinite M\"obius band (figure by C. Varea). Regarded as a  vector bundle 
over $\mathbb R\mathbb P^1$, the space $ \widetilde{\mathbb R^2}$ is used 
within the axioms of characteristic classes 
of real vector bundles for the nontriviality axiom of Stiefel--Whitney classes.

Presented axiomatically the theory of Chern classes looks rather abstract. 
A more geometric approach is to write 
Chern classes in terms of  polynomials on the curvature of a connection, as follows.
Choosing any connection $A$ on a complex vector bundle $E$, the first Chern 
class is actually the cohomology class of the curvature of the connection, that is:
$$c_1(E) = \left[\frac{\sqrt{-1}}{2\pi}F_A\right]  \in H^2(M)$$ 
and the higher Chern class $c_i$  is given by the class of 
the elementary invariant polynomial $P^i$ of degree $i$
$$c_i(E) = \left[P^i\left(\frac{\sqrt{-1}}{2\pi}F_A\right)\right]  \in H^{2i}(M),$$ see \cite[p.\thinspace 407]{GH}.

The Gauss--Bonnet theorem then tells us that  Chern classes measure obstruction to constructing nowhere vanishing sections.
For example, the top Chern class gives the number of zeros of a generic section. In further generality, the 
$i$-th Chern class of a rank $r$ vector bundle is Poincar\'e dual to the degeneracy cycle of $r-i+1$ sections \cite[p.\thinspace 413]{GH}.

\section{Local surfaces and local characteristic classes}

We will consider the local surfaces  $Z_k\ce \Tot \mathcal O_{\mathbb P^1}(-k)$. We begin by setting some notation:

We fix once and for all coordinate charts on $Z_k$, which we refer to as {\it canonical coordinates}, given by
\begin{equation*}
U = \mathbb C^2_{z,u} = \bigl\{(z,u)\bigr\} \qquad\text{and}\qquad
V = \mathbb C^2_{\xi,v} = \bigl\{(\xi, v)\bigr\} \text{,}
\end{equation*}
such that on $U \cap V = \mathbb C^* \times \mathbb C$ we identify
\begin{align}
\label{identification}
\boxed{(\xi, v) = (z^{-1}, z^k u)}\text{.}
\end{align}

We denote by $\ell$ the subvariety of $Z_k$ corresponding to the zero section of 
$\mathcal O(-k)$, thus  $\ell\simeq \mathbb P^1$. Let  $X_k$ denote 
the surface obtained from $Z_k$ by  contracting $\ell$  to a point $x$, with 
$$\pi\colon Z_k \rightarrow X_k$$ the contraction map. Then, $X_1\simeq \mathbb C^2$ and 
$X_k$ is singular at $x$ for $k>1$.

If $E$ is a vector bundle on $Z_k$ then the direct image 
$\pi_*E$ happens to be  a vector bundle on $X_k$ only  in the case when 
$E\vert_\ell $ is trivial, otherwise, it is a non locally free sheaf with the stalk at $x$ strictly larger than at 
the smooth points of $X_k$.

We will consider rank 2 vector  bundles and we wish to define the  local second Chern class of a bundle $E$
around $\ell$ using information from the sheaf $\pi_*E$. This local second Chern class will correspond to 
the local charge of an instanton on $Z_k$.

Considering   $\pi_*E$ then takes us to the realm of 
Chern classes of sheaves on singular varieties. 
Quoting T. Suwa we observe that 
for a  singular complex algebraic or analytic variety $X$, there
are at least three (in general different) kinds of Chern classes in the cohomology $H^*(X)$.
\begin{itemize}
\item the Chern--Schwartz--MacPherson class $c_* (X)$, 
\item the 
Chern--Mather
class $c^M (X)$ and 
\item the canonical class or 
Fulton--Johnson's Chern class
$c^{FJ}(X)$. 
\end{itemize}
These three classes all reduce to 
the standard Chern classes $c(X)$ when the variety
has no singularities. Since there is a priori no reason to prefer one of these definitions over the others,
 we chose instead
to consider the concept of {\it local holomorphic Euler characteristic} $\chi(\ell, E)$ whose definition 
is uncontroversial. 

If $\pi\colon (\widetilde{X}, \ell) \rightarrow (X,x)$ is a resolution of an isolated singularity, 
and $\mathcal F$ is a sheaf of rank $n$ on $\widetilde{X}$, then the {\it local holomorphic Euler 
characteristic} of $\mathcal F$ is 
$$\chi(x,\mathcal F) = \chi(\ell, \mathcal F) = h^0(x,Q) + \sum_{i=1}^n (-1)^ih^0(X, R^i\pi_*(\mathcal F)) $$
where $h^0$ is the dimension of the $0$-th \v{C}ech cohomology,
 $R^i\pi_*(\mathcal F)$ is the $i$-th higher derived image of $\mathcal F$ (see \cite[p.\thinspace 250]{H}) and $Q$ is 
the skyscraper sheaf supported at $x$ defined by the exact sequence
$$0 \rightarrow \mathcal \pi_*F \rightarrow (\pi_*\mathcal F)^{\vee\vee} \rightarrow Q \rightarrow 0.$$

 If $X$ is  a compact orbifold and $\mathcal F$ a sheaf over $X$, there exists 
the global  {\it holomorphic orbifold Euler characteristic} of $\mathcal F$  (see \cite{Bl})
$$\chi_{orb}(X, \mathcal F) = \int_X  \mbox{ch}(\mathcal F) \mbox{td}(X).$$
It satisfies:
$$\chi(X, \mathcal F)= \chi_{orb}(X, \mathcal F) + \sum_{x \in \Sing(X)}\mu(x,\mathcal F).$$ 
Here $\mu(x,\mathcal F)$ is a rational number that depends both on the class of the sheaf 
$\mathcal F$ in the Grothendieck 
group of the variety and on the order of singularity, 
with the denominator at an orbifold singularity  given by the order of the corresponding group
that gives the local quotient.

In the examples we wish to consider when $E $ is  a rank 2 vector bundle on $Z_k$ then we refer to 
the local holomorphic Euler characteristic of $E$ around $\ell$ as the {\it charge} of $E$ for short, 
since it corresponds to the contribution to the  charge in the language of instantons. It is given by 
$$\chi(\ell, E ) = h^0(Q) + h^0(R^1\pi_*E)$$
where 
$$0 \rightarrow \pi_*E \rightarrow (\pi_*E)^{\vee\vee} \rightarrow Q \rightarrow 0.$$
The  numbers $h^0(Q)$ and  $h^0(R^1\pi_*E)$ are independent analytic invariants, called the {\it width}
and the {\it height} of $E$ respectively.
These 2 invariants can be calculated using the theorem on formal functions (see \cite[p.\thinspace 277]{H})
$$R^if_*(E)^{\wedge}_x= \lim_{\leftarrow}H^i(\ell^n, E_n)$$
by calculating \v{C}ech cohomology in canonical coordinates as given in (\ref{identification}).
It is the local holomorphic characteristic $\chi(\ell, E)$ 
that gets identified with the local instanton charge under the Kobayashi--Hitchin correspondence,
and for this reason we call it the {\it charge} of the bundle, it gives the local contribution to the second Chern class
$c_2^{loc}(X,E)$, which we will denote simply by $c_2^{loc}(E)$ when the base $X$ is clear from the context.

To provide examples of these invariants, we first recall some properties of vector bundles on local surfaces.
To start with, holomorphic vector bundles on $Z_k$ are filtrable and algebraic
\cite[Lem.\thinspace 3.1, Thm.\thinspace 3.2]{G1}. Filtrability for a rank 2 bundle $E$ means that it is an 
extension of line bundles. While filtrability always happens for bundles on complex curves (Riemann surfaces) 
this is a very unusual property for complex surfaces. For example, on the projective space $\mathbb P^2$ only 
split bundles are filtrable, so that the entire moduli spaces of stable 
vector bundles on $\mathbb P^2$ are made
of nonfiltrable bundles.

 If $E$ is a rank 2 bundle on $Z_k$ with $c_1(E)=0$ then for some integer $j\geq 0$
 there is a short exact sequence
\begin{equation} \label{ext} 0 \rightarrow 	\mathcal O(-j)\rightarrow E \rightarrow \mathcal O(j) \rightarrow 0,\end{equation}
where $E\vert_\ell = \mathcal O(-j) \oplus \mathcal O(j)$ splits by Grothendieck's splitting principle.
The integer $j$ is called the \emph{splitting type} of $E$. Note that here we use the symbol
$\mathcal O(j)$ for the line bundle with first Chern class $j$ both on $\ell\simeq \mathbb P^1$ 
and on $Z_k$, but there should be no confusion given that one is obtained from the other by pullback.

Giving a short exact sequence $0 \rightarrow 	\mathcal O(-j)\rightarrow E \rightarrow \mathcal O(j) \rightarrow 0$
corresponds to giving an extension class
 $p \in \mathrm{Ext}^1(\mathcal O(j), \mathcal O(-j))$. This extension class can be 
 written in canonical coordinates as a polynomial such that the transition matrix of $E$ in canonical coordinates is 
	$\left(\begin{matrix} 
		z^j & p \\
        0 & z^{-j}
	\end{matrix}\right).$
Therefore,  the bundle $E$ is completely determined by the data $(j,p)$.    

The following table contains examples for the values of width, height and charge for vector bundles on
 $Z_1$ with splitting type $j=3$ defined by the monomials in the first column:
\begin{table}[ht]	
\begin{center}
\begingroup
\renewcommand*{\arraystretch}{1.15}
\begin{tabular}{c||c|c|c}
monomial  & {width} &{height} &{charge} \\
\hline
$z^{-1}u$ & 3 & 2 & 5\\
$u$ & 1 & 2 & 3\\
$zu$ & 1 & 2 & 3\\
$z^{2}u$ & 3 & 2 & 5\\
\hline
$u^2$ & 3 & 3 & 6\\
$zu^2$ & 2 & 3 & 5\\
$z^{2}u^2$ & 3 & 3 & 6\\
\hline
$zu^3$ & 4 & 3 & 7\\
$z^{2}u^3$ & 4 & 3 & 7\\
\hline
$z^{2}u^4$ & 5 & 3 & 8\\
\hline
zero & 6 &3 &  9\\
\end{tabular}
\endgroup
\vspace{2mm}
\end{center}
\caption{$j = 3$ on $Z_1$}
\label{table}
\end{table}

\section{Moduli spaces}

To define moduli of vector bundles, we first fix topological invariants, and then mod out by 
holomorphic isomorphisms. Next, as  it is well known in algebraic geometry,
 we need to  select a good subset of bundles to
obtain moduli spaces. This means that we have   to select a proper subset 
of the set of all vector bundles with the chosen topology to obtain a well behaved quotient,
 else the quotient may be a stack, or something worse, but will not be an algebraic 
 variety, at times not even a complex space. 
In geometric invariant theory, a well behaved moduli space is obtained using the 
concept of stable bundle. However, the concept of stability is defined for bundles on compact 
varieties and works well in such case, but it does not apply to bundles on
 noncompact varieties. For this reason, a somewhat ad-hoc choice of 
 good moduli space was obtained for the surfaces $Z_k$ by choosing those bundles
 that do not split on the first formal neighborhood of $\ell$. In practice this means 
 that the extension class has a nontrivial linear term on $u$, where $u=0$ 
 is the equation that cuts out $\ell$  inside the surface. Once this is done, 
 then \cite[Thm.\thinspace 4.11]{BGK} obtained: 
 The {\it moduli space of holomorphic bundles} of rank 2  on $Z_k$ 
with $c_1=0$ and splitting type $j\geq k\geq 2$ is a quasi-projective variety of dimension $2j-k-2$.

Gasparim, K\"oppe, and Majumdar \cite[Cor.\thinspace 5.5]{GKM} showed that 
a rank 2 bundle on $Z_k$ corresponds to an $\mathrm{SU}(2)$ instanton if and only if its splitting type is a multiple of $k$. 
The correspondence also takes gauge transformation of instantons to isomorphism of holomorphic bundles. 
Furthermore, their choice of ad-hoc stability concept is such that {\it stable} instantons with splitting type $j$ have 
topological charge $j$. 
It then follows that the {\it moduli space of  $\mathrm{SU}(2)$ instantons}  of charge $j=nk\geq 2$ on $Z_k$ is a 
quasi-projective variety of dimension $2nk-k-2$. There is a single instanton of charge $1$ on $Z_1$ and, 
very surprisingly,  there are no charge $1$ instantons on $Z_k$ for $k\geq 2$, as we shall see next. 

%
%
%\begin{theorem}[Ballico, G., K\"oppe]
%The invariants width and height stratify the moduli stack of bundles on $Z_k$ into Hausdorff components.
%\end{theorem}
%

\section{Instanton decay}

Sharp bounds for local charges of $\mathrm{SU}(2)$  instantons on $Z_1$ with splitting type $j$ are as follows:
$$ j \leq c_2^{loc}(Z_1, E)\leq j^2.$$
The minimum local charge of a splitting type $j$ instanton on $Z_1$ equals $j$ and is attained by those 
bundles that do not split on the first formal neighborhood, such as those with local data $(j, p=a u + b zu)$ 
for any constants $a,b $ not both $0$.
The maximum local charge is always attained by the split bundle, thus the case of the  instanton with data $(j, p=0)$ 
which has local charge $j^2$.
Furthermore, every pair of admissible numerical values of $(width, height)$ occurs for some  instanton  on $Z_1$ \cite[Thm.\thinspace 0.2]{BG}.
Therefore, there exist instantons on $Z_1$ for any values of local charge, width and height. 

In contrast, existence of instantons on  $Z_k$ with prescribed value of local charge is still an  open problem.
Indeed,  even 
though sharp bounds are know for each fixed splitting type, it is not know whether all intermediate values 
of width and height actually occur. This issue is already featured in algebraic geometry, since the corresponding 
existence of vector bundles is unknown. The current state of affairs is as follows. 
Assume  $k\geq 2$. Let $E$ be a  holomorphic rank 2 bundle on $Z_k$  with  $c_1=0$ and splitting type $j= nk+r$
with $0 \leq r < k$, then the following bounds are sharp \cite[Cor.\thinspace 2.18]{BGK}:
\begin{equation}\label{bounds} j-1 \leq c_2^{loc}(Z_k, E) \leq \left\{\begin{array} {lll} n^2k + r(2n+1)-1 & if & r \neq 0 \\
                                                                              n^2k & if & r=0.\end{array} \right. \end{equation}

Even though the question of existence of instantons on $Z_k$ has not been completely  understood when $k\geq 2$,
the current partial knowledge generated strong conclusions regarding the possibilities for instanton decay. 
The main result obtained in \cite[Thm.\thinspace 6.8]{BGK} is that if 
 $k\geq 2$, then 
the minimal local charge of a nontrivial instanton on $Z_k$ is $k-1$.
As a corollary we conclude that  instanton decay is obstructed by the self-intersection $\ell^2=-k$ 
of the line in $Z_k$, given that there are no instantons of charge 1 on $Z_k$ when $k>2$. \\

\noindent{\bf Open question}: For what admissible values of local charge  $c$ fitting into bounds (\ref{bounds})
do there exist instantons on $Z_k$ of charge $c$?\\

\noindent{\bf Open question}: Assume $A$ is an $\mathrm{SU}(2)$ instanton on $Z_k$ of minimal charge, hence the local charge 
satisfies $c(A)= k-1>1 $ if $k>2$. Then, such an instanton $A$ ought to decay, because maintaining
such local charge costs too much energy. However, $A$ can not decay to other instantons on $Z_k$, because 
there are no nontrivial instantons of lower charge in this space. What is the result of the decay 
of such an instanton? We will give a partial solution to this question in the following section.

\section{Holomorphic surgery and instanton decay}

We define a new type of holomorphic surgery which we call {\it grafting}. 

\begin{center}
\begin{tikzpicture}

\node (P1) at (-2.3,2) {\includegraphics[scale=.17]{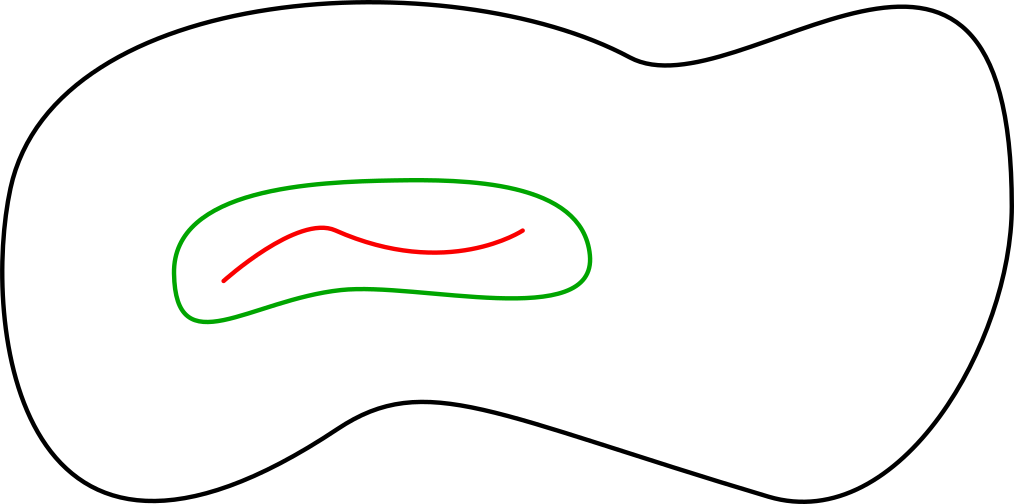}};

\node (P2) at (4,2) {\includegraphics[scale=.17]{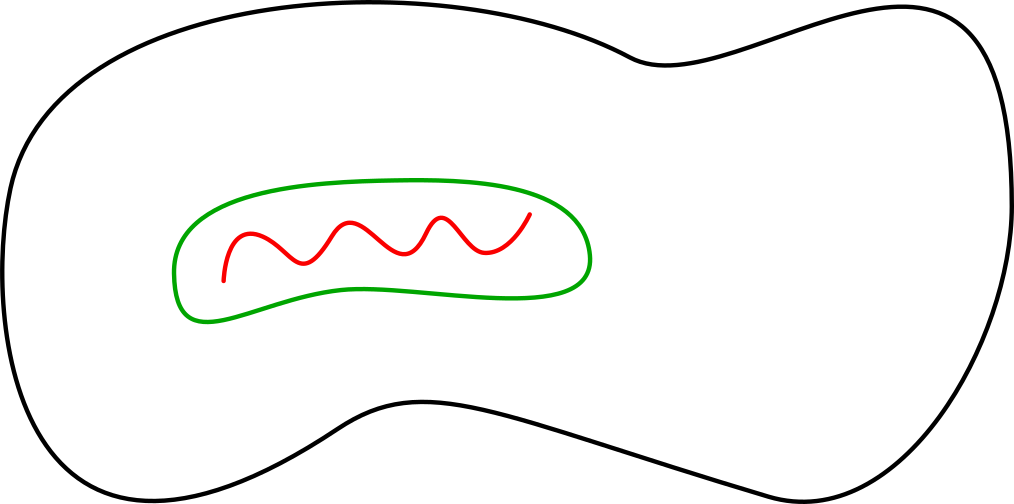}};

\node (P3) at (-2.3,-2) {\includegraphics[scale=.17]{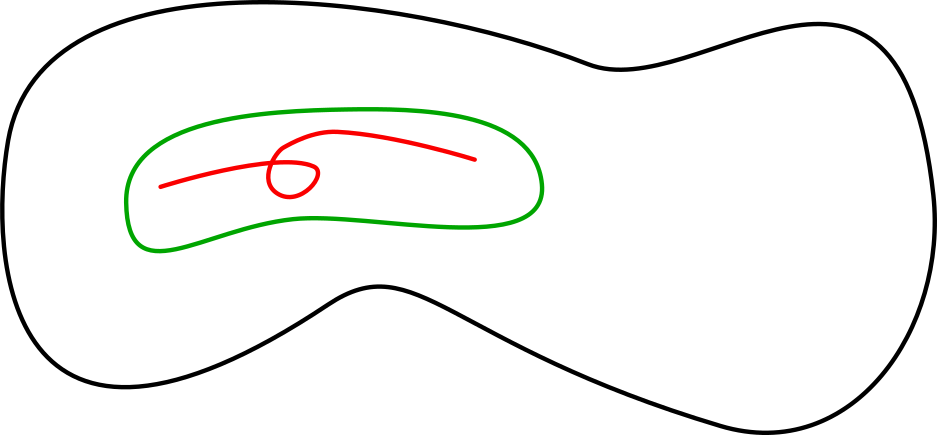}};

\node (P4) at (4,-2) {\includegraphics[scale=.17]{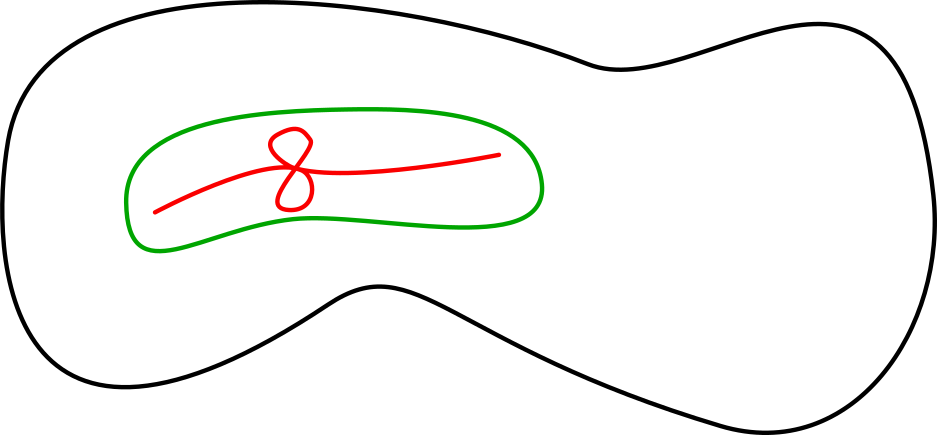}};

\draw[thick, ->] (P1) -- (P3);
\draw[thick, ->] (P2) -- (P4);
\draw[thick, ->] (P1) -- (P2) node[midway,above] {\textcolor{OliveGreen}{graft}};

\node[OliveGreen] (Q1) at (-1.4,2) {$Z_k$};
\node[OliveGreen] (Q2) at (4.7,2) {$Z_{k'}$};
\node[OliveGreen] (Q3) at (-1.4,-2) {$X_k$};
\node[OliveGreen] (Q4) at (4.7,-2) {$X_{k'}$};

\end{tikzpicture}
\end{center}

We will discuss 2 types of grafting: open and closed. 
Closed grafting is a specific new type of surgery on the base space such that, when the space holds an instanton (or a vector bundle),
then the 4 dimensional manifold changes under grafting, but the instanton (or holomorphic bundle) on it remains essentially 
the same, although loosing some of its charge. 
We will carry out the closed grafting by   replacing a $Z_k$ with  another $Z_{k'}$. 
Holomorphic patching of a $Z_1$ is possible at any smooth point, and corresponds to the operation of 
blowing up a point which is well known in algebraic geometry.
 However, the new operation considered here is  not just that of changing the surface by adding a $Z_1$,
but also carrying along with it an instanton (or holomorphic vector bundle). 

We may regard the blow up  $\pi \colon \widetilde X \rightarrow  X $ of a point $x \in X$ as a holomorphic patching
$$\widetilde{X} = (X -\{x\})  \cup_{ Z_1^o} Z_1.$$
Given holomorphic bundles  $E$  on $X$ and $F$  on $Z_1$ with the same rank,  we may
construct an open grafting of $E$ and $F$, written $E \atrevo F$, read ``$E$ graft $F$", by making
$E \atrevo F$ isomorphic to $E$ on $\widetilde{X}- \ell$ and isomorphic to $F$ on $Z_1$;
this requires  specifying a choice of gluing. So, the holomorphic type of 
the resulting bundle depends on the  gluing,  but the resulting second Chern class of $E \atrevo F$ depends only 
on $E$ and $F$, see \cite{G2}. 
It is worth noticing that such grafting is always possible, because of the fact that 
every holomorphic bundle on $Z_1$ is  trivial on $Z_1^o\ce Z_1 - \{\ell\}$, see \cite{G1}. 
Similarly, we can graft a bundle locally by exchanging it on a  $Z_k$  type of neighborhood inside a complex surface.
We  have 
$$c_2(E \atrevo F) = c_2(E) + c_2^{loc}(F).$$

The second type of grafting, which is a somewhat more surprising operation is to essentially keep the vector bundle while
changing the base manifold. Suppose that we have a complex surface $X$ that contains a $-k$ line. Hence, 
$X$ contains a line $\ell$ with a tubular neighborhood  of the same holomorphic type as a neighborhood of $\ell$ in $Z_k$. 
We now replace $k$ by $k'$ by cutting out the $-k$ line and gluing in a $-k'$ line. If $E$ is a bundle on $X$, we keep all of its 
transition matrices, and they now define a new bundle on $(X -\ell)  \cup_{Z_k^o} Z_{k'}$
which we denote by $E \trevo _{k'}$ read ``$E$ graft $k'$".  In general we will need to also modify other charts of $X$ because 
of the cocycle condition in triple intersections. But, let us consider the simplest case 
when $X$ is one of the $Z_k$ themselves, so that we just exchange $Z_k$ by $Z_{k'}$ and 
keep the transition matrix for the bundle.  This might 
seem a trivial thing to do,  nevertheless, the effect on the local charge is very strong, as the following examples illustrate. 

First observe  that over $Z_k$ the bundle 
$E=\mathcal O(j) \oplus \mathcal O(-j)$, that is, the split bundle with type
 $j=nk$ and $p=0$ over $Z_k$, has 
local charge $c_2^{loc}(Z_k, E) = n^2k$, in fact, it realizes the upper bound in \ref{bounds}. Now, let us 
see the effect of grafting in a couple of cases,
 by considering the split bundles with type $j=6$ over $Z_1$ and $Z_2$.
So that, in each case, we are considering a bundle defined by transition matrix 
  $T= \left(\begin{matrix} z^6 & 0 \\
                                    0 & z^{-6}
           \end{matrix}\right).$

\noindent{\sc Case 1}. When we consider the  split bundle $E\rightarrow Z_1$  given by 
transition matrix $T$  we are in the case $j=6\times 1$, so 
 $$c_2^{loc}(E)= 6^2 \times 1= 36$$ by grafting to $Z_3$ we get the bundle 
 $E\trevo_3$ also given by transition matrix $T$ but now we are in the case
 $j= 2\times 3$, so 
 $$c_2^{loc}(E\trevo_3)= 2^2 \times 3= 12.$$
 So that in this grafting the instanton bundles lost 24 units of charge.\\
 
\noindent{\sc Case 2}. When we consider the  split bundle $E\rightarrow Z_2$  given by 
transition matrix $T$  we are in the case $j=3\times 2$, so 
 $$c_2^{loc}(E)= 3^2 \times 2= 18$$ by grafting to $Z_6$ we get the bundle 
 $E\trevo_6$ also given by transition matrix $T$ but now we are in the case
 $j= 1\times 6$, so 
 $$c_2^{loc}(E\trevo_6)= 1^2 \times 6= 6.$$
 So that in this grafting the instanton bundles lost 12 units of charge.\\

In both cases we can interpret the loss of charge as the energy used to 
provoke the grafting in the base. 

Now let $j \geq 2$ and $E \to Z_1$ the bundle defined by the transition matrix $T' = \left( \begin{matrix}
z^j & zu \\ 0 & z^{-j}
\end{matrix} \right)$. As $E$ does not split, we have that 
\[
c_2^{loc}(E) = j.
\]
By grafting to $Z_3$ we have that
\[
c_2^{loc}(E \trevo_3) = j-1,
\]
so in this grafting the instanton bundles lost only 1 unit of charge.

Here is  a picture of closed grafting, as it exchanges the 
local surfaces $Z_k$ and $Z_{k'}$ inside a complex surface,
while essentially maintaining the instanton bundle  
globally. The grafted bundle  is still given by the same transition matrices, however, loses global charge, 
that is, undergoes a loss of  second Chern class.

\begin{center}
\begin{tikzpicture}

\node (P1) at (-2,2) {\includegraphics[scale=.17]{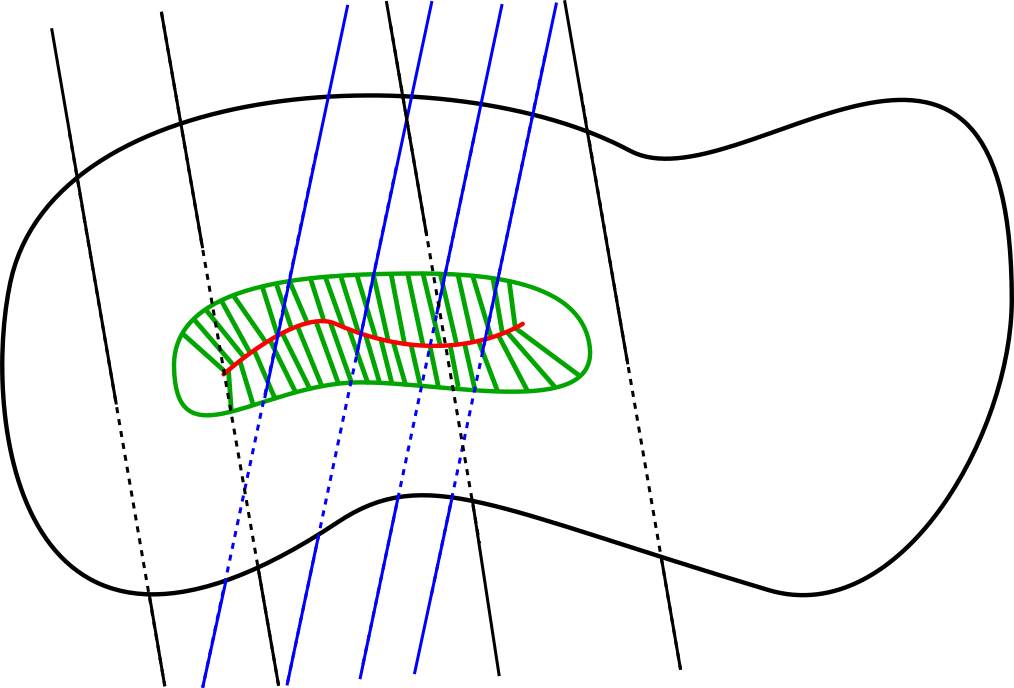}};

\node (P2) at (3.8,2) {\includegraphics[scale=.17]{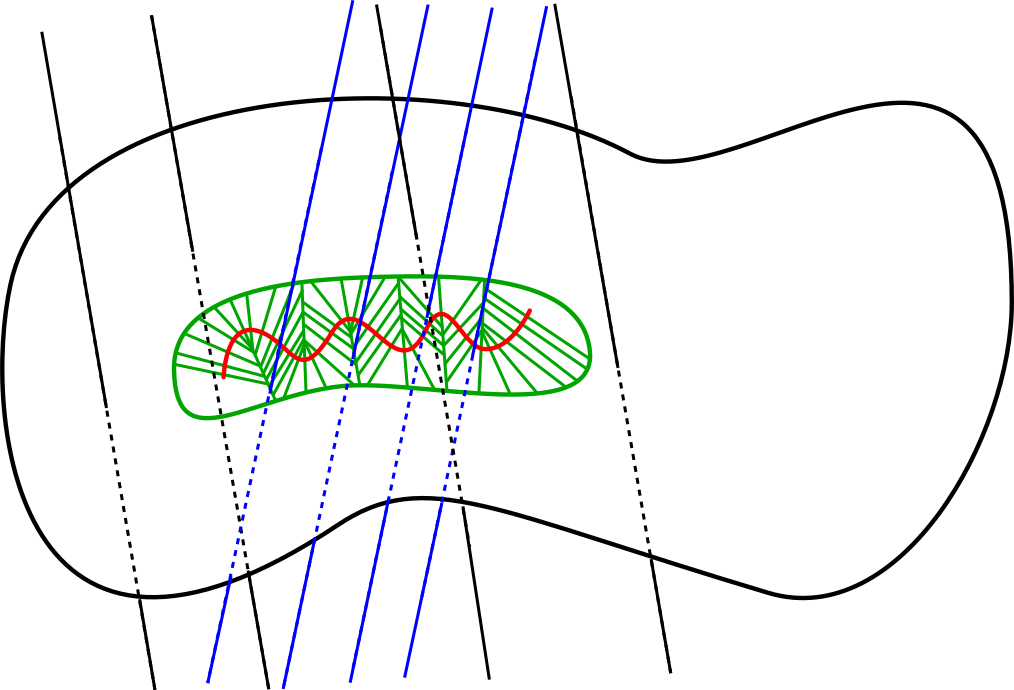}};

\node (P3) at (-2,-2) {\includegraphics[scale=.17]{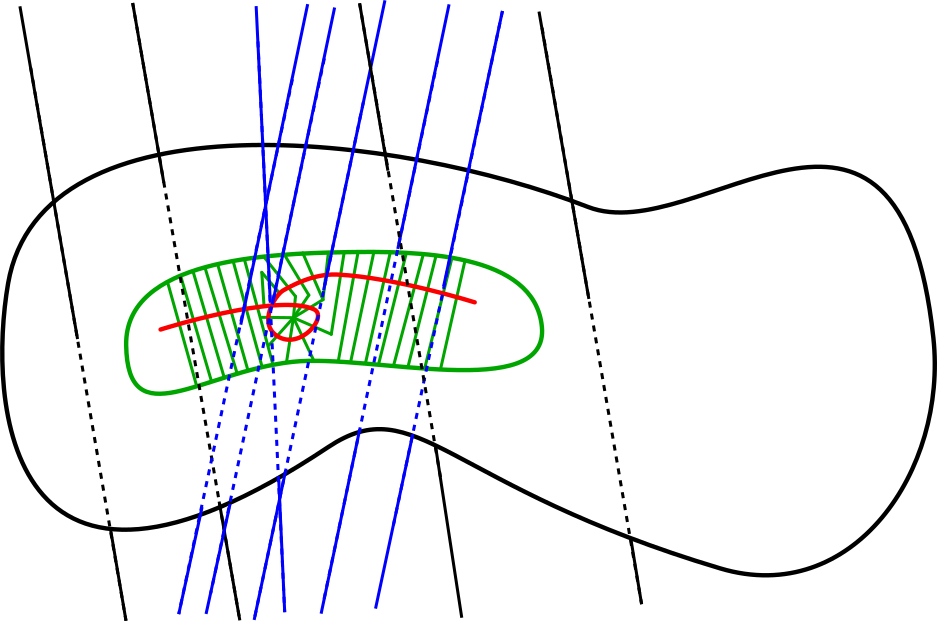}};

\node (P4) at (3.8,-2) {\includegraphics[scale=.17]{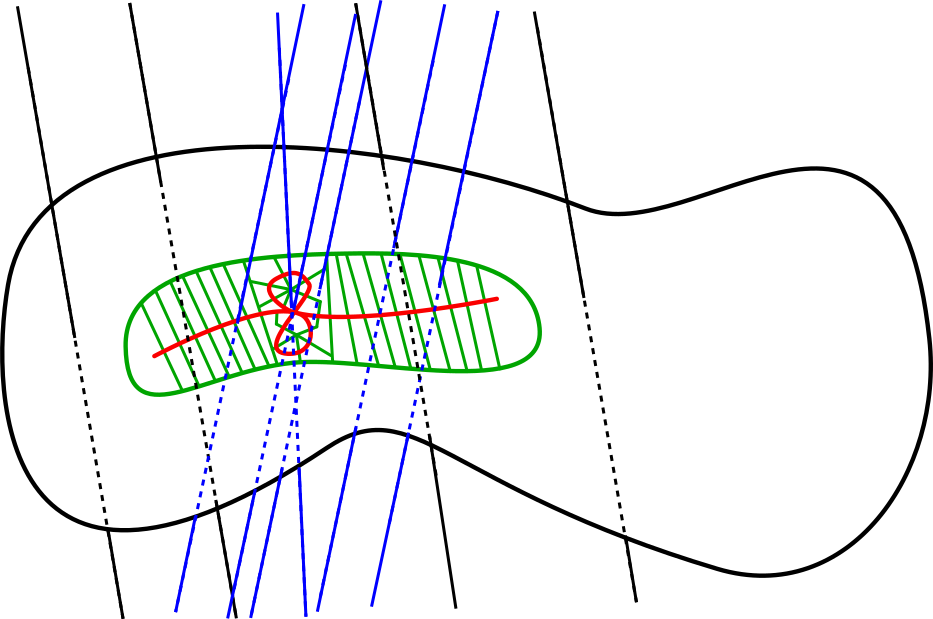}};

\draw[thick, ->] (P1) -- (P3);
\draw[thick, ->] (P2) -- (P4);
\draw[thick, ->] (P1) -- (P2) node[midway,above] {\textcolor{OliveGreen}{ graft}};

\node[OliveGreen] (Q1) at (-1.1,2) {$Z_k$};
\node[OliveGreen] (Q2) at (4.7,2) {$Z_{k'}$};
\node[OliveGreen] (Q3) at (-1.1,-2) {$X_k$};
\node[OliveGreen] (Q4) at (4.7,-2) {$X_{k'}$};

\end{tikzpicture}
\end{center}

In conclusion, we propose that an instanton with high charge provokes a 
grafting on the 4 dimensional base manifold which holds it. The grafting happens around a 
2-sphere where the instanton charge is highly concentrated.  The instanton looses 
charge in the process, while the curvature of the base manifold around the 2-sphere increases. 
In other words, we have explained that the decay of an instanton  can 
happen when the instanton inflicts larger  concentrations of curvature 
 around 2-spheres, thereby losing charge. We then say that the addition to the curvature is 
 grafted by the instanton.

\section{Acknowledgment}

 I first met Pushan Majumdar in 1998 and we immediately engaged into scientific discussions. 
 Shortly after that,  I found him reading about moduli of instantons, 
and inquired why he was reading about such a theme, which was not part of his PhD project. He said 
that he wanted to be my friend, and if he did not learn about the themes that interested me, 
then soon we would not have anything to talk about. Decades later, this remains the most loyal offer 
of friendship I have ever received. E.G. 

We thank Koushik Ray for inviting us to contribute to the memorial volume of Pushan and for 
suggesting several improvements to our article.

\end{document}